\documentclass[12pt]{article}
\usepackage{amssymb,latexsym}
\input xypic \xyoption{all}
\title{Nahm transform for doubly-periodic instantons}
\author{Marcos Jardim \\ University of Pennsylvania \\ Department of Mathematics \\
Philadelphia, PA 19104-6395 USA \\ jardim@math.upenn.edu} 

\newcommand{\pf}{{\em Proof: }} \newcommand{\pfend}{\hfill $\Box$ \linebreak}
\newcommand{\eq}{\begin{equation}} \newcommand{\seta}{\rightarrow}
\newcommand{\torus}{T\times\cpx} \newcommand{\tproj}{T\times\proj}
 
\newcommand{\dual}{\hat{T}} 
\newcommand{\del}{\overline{\partial}} 
\newcommand{\see}{\Leftrightarrow} \newcommand{\imply}{\Rightarrow}
 \newcommand{\ksi}{\xi} 
\newcommand{\real}{\mathbb{R}} \newcommand{\cpx}{\mathbb{C}}

\newcommand{\proj}{\mathbb{P}^1} \newcommand{\vv}{{\cal V}} \newcommand{\ee}{{\cal E}}

\newcommand{\cala}{{\cal A}} \newcommand{\calg}{{\cal G}} \newcommand{\jj}{{\cal J}}
\newcommand{\ind}{{\rm index}}  
 \newcommand{\as}{\pm\ksi_0}
\newcommand{\oo}{{\cal O}}  
\newcommand{\cfg}{{\cal A}_{(k,\ksi_0)}} \newcommand{\matriz}{\left( \begin{array}{cc}}
\newcommand{\column}{\left( \begin{array}{c}} \newcommand{\matfim}{\end{array} \right)}

\newtheorem{thm}{Theorem} \newtheorem{lem}[thm]{Lemma} \newtheorem{prop}[thm]{Proposition}

\begin{document}
\maketitle

\begin{abstract}
We present the Nahm transform of the doubly-periodic instantons
previously introduced by the author, converting them into certain 
meromorphic solutions of Hitchin's equations over an elliptic curve.
\end{abstract}

\baselineskip18pt \newpage

\section{Introduction} \label{intro}

This paper is the second in a series of three. In the first paper
\cite{J2}, we have shown how $SU(2)$ certain instantons over
$\real^4$ which are periodic in two directions, so called {\em
doubly-periodic instantons}, can be constructed from a particular
type of singular solutions of Hitchin's equations over an elliptic
curve \cite{H}. This was done via a procedure known as {\em Nahm 
transform}, which has attracted much attention among physicists 
recently (see for instance \cite{KS}, \cite{vB} and the references 
therein).

We now present the inverse construction, showing that all extensible
doubly-periodic instantons were obtained in \cite{J2}.

Recall that given a function $f:\cpx\seta\real$, we say that
$f\sim O(|w|^n)$ if:
$$ \lim_{w\seta\infty} \frac{|f(w)|}{|w|^n} < \infty $$

We consider anti-self-dual connections $A$ on rank two bundle
$E\seta\torus$ satisfying the following conditions:
\begin{enumerate}
\item $|F_A|\sim O(r^{-2})$;
\item there is a holomorphic vector bundle $\ee\seta\tproj$ with trivial
determinant such that $\ee|_{T\times(\proj\setminus\{\infty\})}\simeq(E,\del_A)$,
where $\del_A$ is the holomorphic structure on $E$ induced by the instanton
connection $A$;
\end{enumerate}
Such connections are said to be {\em extensible}. Moreover, we assume
the restriction of the extended bundle to the added divisor splits as a sum
of flat line bundles, i.e.:
$$ \ee|_{T_\infty}=L_{\ksi_0}\oplus L_{-\ksi_0} $$
and $\pm\ksi_0$ can be seen as points in the Jacobian torus
$\jj(T)=\dual$. We say $\ksi_0$ is the {\em asymptotic state} of
the connection $A$. We also fix the topological type of the
extended bundle $\ee$ by making $c_2(\ee)=k>0$; the integer $k$ is
the {\em instanton number} of the connection $A$. The space of such
connections is denoted by $\cala_{(k,\xi_0)}$.

Let us now outline the contents of this paper. The key feature of
Nahm transforms is to try to solve the Dirac equation, and then use
its $L^2$-solutions to form a vector bundle over the jacobian torus
$\dual$, which parametrises the set of holomorphic flat line
bundles over $\torus$. Therefore, our first task is to show that
the Dirac operator is Fredholm and compute its index.

The bulk of the paper lies in sections \ref{nahm} and \ref{holo},
where we present the Nahm transform of doubly-periodic instantons
and show some of the properties of the transformed objects.

Section \ref{inv} is dedicated to prove that the construction here
presented is indeed the inverse of the one presented in \cite{J2},
completing the proof of the main result that motivated this series
of papers:

\begin{thm} \label{nahmthm}
The Nahm transform is a bijective correspondence between the following objects:
\begin{itemize}
\item gauge equivalence classes of extensible, irreducible $SU(2)$
instanton connections on $E\seta\torus$ with fixed instanton number
$k$ and nontrivial asymptotic state $\ksi_0$; and
\item admissible $U(k)$ solutions of the Hitchin's equations over
the dual torus $\dual$, such that the Higgs field has at most
simple poles at $\as\in\dual$, with semi-simple residues of rank $\leq2$
if $\ksi_0$ is an element of order 2 in the Jacobian of $T$, and
rank $\leq1$ otherwise.
\end{itemize} \end{thm}

We conclude this article by stating a higher rank generalization of
the above result.

In the third and last paper in this series \cite{J3}, we shall
discuss the role played by {\em spectral curves} in the
correspondence between doubly-periodic instantons and singular
solutions of Hitchin's equations, thus completing a circle of ideas
analogous to Hitchin's approach to monopoles \cite{H3}:
\vskip18pt 
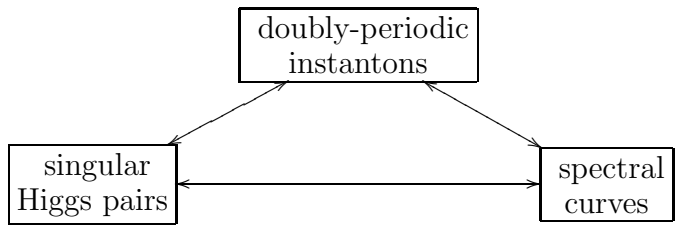
\centerline{
\xymatrix{
& *+[F]\txt{ doubly-periodic\\ instantons }
\ar@{<->}[dl] \ar@{<->}[dr] & \\
*+[F]\txt{ singular \\ Higgs pairs }
\ar@{<->}[rr] & &
*+[F]\txt{ spectral \\ curves }
}}
\vskip18pt

\paragraph{Acknowledgements.}
This work is part of my Ph.D. project \cite{J}, which was funded by
CNPq, Brazil. I am grateful to my supervisors, Simon Donaldson and
Nigel Hitchin, for their constant support and guidance. I also thank
Brian Steer and Olivier Biquard for valuable suggestions in the
later stages of this project.

%-------------------------------------------------------------------

\section{Extensibility and asymptotic behaviour}

We now use the extensibility hypothesis to study the compatibility
between the instanton connection $A$ and the extended bundle
$\ee\seta\tproj$. More precisely, we first want to show that the
holomorphic type of the restriction of the extended bundle to the
added divisor  $T_\infty=T\times\{\infty\}$ is indeed directly
determined by the asymptotic behaviour of the instanton connection
$A$. Then we argue that the topology of $\ee$ is fixed by the
energy ($L^2$-norm) of $A$.

Before that, we must fix an appropriate trivialisation at infinity.

\paragraph{Gauge fixing at infinity.}
Let $B_R$ denote a closed ball in $\cpx$ of radius R, and let $V_R$
be its complement. Also, consider the obvious projection 
$p:T\times V_R\seta T$. We shall need the following technical proposition, 
which follows from the gauge-fixing result established in \cite{BJ} 
(see also the appendix in \cite{J}).

\begin{prop} \label{goodgauge}
If $|F_A|\sim O(r^{-2})$, then, for $R$ sufficiently large, there
is a gauge over $T\times V_R$ and a constant flat connection $\Gamma$
on a topologically trivial rank two bundle over the elliptic curve
such that:
$$ A-p^*\Gamma \ = \ \alpha \ \sim O(r^{-1}\cdot\log r) $$
\end{prop}

%----------------------------------------------------

\subsection{Asymptotic states}

By general theory, a constant flat connection
on a bundle $S\seta T$ determines uniquely a holomorphic structure on this
bundle. Moreover, $S$ must split, holomorphically, as the sum of two line
bundles, i.e. $S=L_{\ksi_0}\oplus L_{-\ksi_0}$, uniquely up to $\pm1$. Here,
$\as$ are seen as points in $\dual$, the Jacobian of the elliptic curve $T$.

Therefore, by proposition \ref{goodgauge}, to each extensible instanton
connection we can associate an unique pair of opposite points $\as\in\dual$.
Such points are called the {\em asymptotic states} of $A$.

\begin{lem}
If an extensible instanton connection $A$ has asymptotic states $\as$, then
$\ee|_{T_\infty}=L_{\ksi_0}\oplus L_{-\ksi_0}$.
\end{lem}

\pf
Let $V_\infty\subset\proj$ be a small neighbourhood centred at
$\infty\in\proj$; let $w$ be a coordinate there. We can regard
$\ee|_{T\times V_\infty}$ as a family of rank 2 bundles over $T$,
parametrised by $w$, Furthermore,
If $\del$ denotes the holomorphic structure on $\ee$, let
$\del_w$ be the holomorphic structure on the restriction $\ee|_{T_w}$.
Clearly, as operators:
$$ \lim_{w\seta\infty}\del_w=\del_\infty $$
However, from condition (2) in the definition of extensibility,
we know that $\del_w=\del_{A|_{T_w}}$ away from $\infty$. But
proposition \ref{goodgauge} tells us that $\del_{A|_{T_w}}$
approaches $\del_{\Gamma}$ as $w\seta\infty$. Therefore,
$\del_\infty=\del_{\Gamma}$, and the lemma follows. \pfend

%----------------------------------------------------

\subsection{The instanton number}
Moreover, as we mentioned before, the topological type of $\ee$
is determined by the energy of the instanton connection:

\begin{lem} \label{c2=i}
$c_2(\ee)=\frac{1}{8\pi^2}\int_{\torus}|F_A|^2$
\end{lem}

\pf Again, let $V$ be a small neighbourhood of $\infty\in\proj$. Let
$\Gamma_{\as}$ be the canonical connection on the bundle
$L_{\ksi_0}\oplus L_{-\ksi_0}$ over an elliptic curve and consider the
projection $p:T\times V\seta T$.

Now consider a connection $A'$ on the extended bundle $\ee$ that coincides
with $p^*\Gamma_{\as}$ on $T\times V$. Therefore
\begin{eqnarray} \label{c2}
c_2(\ee) & = & \frac{1}{8\pi^2} \int_{\tproj} {\rm Tr}(F_{A'}\wedge F_{A'}) \ = \
\frac{1}{8\pi^2} \int_{T\times(\proj\setminus\{\infty\})} {\rm Tr}(F_{A'}\wedge F_{A'}) \nonumber \\
& = & \frac{1}{8\pi^2} \lim_{R\seta\infty}\int_{T\times B_R} {\rm Tr}(F_{A'}\wedge F_{A'})
\end{eqnarray}

On the other hand, we have from lemma \ref{goodgauge} that $A-A'=\alpha$ is a
1-form in $O(r^{-1}\cdot\log(r))$. Define the 1-parameter family of connections
$A_t=A'+t\cdot\alpha$, so that the corresponding curvatures:
\begin{eqnarray}
& F_{A_t}=t\cdot F_{A} + (1-t)\cdot F_{A'} -
\left( t-\frac{t^2}{2}\right)\cdot\alpha\wedge\alpha & \nonumber \\
& \Longrightarrow\ \ \ |F_{A_t}|\sim O(r^{-2}\cdot \log^2r)\ \
\forall t\in[0,1] \label{famcurvs}
\end{eqnarray}

So let:
\begin{equation} \label{i}
i(A) = \frac{1}{8\pi^2} \int_{\torus} {\rm Tr}(F_{A}\wedge F_{A}) \ = \
\frac{1}{8\pi^2} \lim_{R\seta\infty}\int_{T\times B_R} {\rm Tr}(F_{A}\wedge F_{A})
\end{equation}
Usual Chern-Weil theory tells us that:
\begin{eqnarray*}
c_2(\ee)-i(A) & = & \frac{1}{8\pi^2} \lim_{R\seta\infty} \left\{ \int_{T\times B_R} \left(
{\rm Tr}(F_{A'}\wedge F_{A'})-{\rm Tr}(F_{A}\wedge F_{A}) \right) \right\} \ = \\
& = & \frac{1}{4\pi^2} \lim_{R\seta\infty} \left\{ \int_{T\times B_R} d \left(
\int_0^1 {\rm Tr}(\alpha\wedge F_{A_t}) \right) \right\} \ = \\
& = & \frac{1}{4\pi^2} \lim_{R\seta\infty} \left\{ \int_{T\times S_R^1} \left(
\int_0^1 {\rm Tr}(\alpha\wedge F_{A_t}) \right) \right\} \ = 0
\end{eqnarray*}
by our estimates in proposition \ref{goodgauge} and in equation 
({\ref{famcurvs}). This completes the proof. \pfend

%----------------------------------------------------------

\subsection{Estimating the Dolbeault operator}
Finally, we need one final lemma that will be useful in the
following section section, where we develop a Fredholm theory for
the Dirac operator coupled to an instanton connection $A\in\cfg$.

First, note that the bundle $L_{\ksi_0}\oplus L_{-\ksi_0}\seta T$ admits a
flat connection with constant coefficients, which we denote by
$\Gamma_{\ksi_0}$. Use the projection $T\times V_R\stackrel{p_1}{\seta}T$ to
pull it back to $T\times V_R$. We show that:

\begin{lem} \label{model}
Let $A\in\cfg$ be any extensible instanton connection. Given
$\epsilon>0$, there is $R$ sufficiently large such that:
$$ ||\del_A-\del_{\Gamma_{\ksi_0}}||_{L^2(T\times V_R)}<\epsilon $$
\end{lem}

\pf Since $\del_A-\del_{\Gamma_{\ksi_0}}$ is just the $(0,1)$-part
of the 1-form $\alpha=A-\Gamma_{\ksi_0}$, the statement is a simple
consequence of the gauge-fixing proposition \ref{goodgauge}. \pfend

%--------------------------------------------------------------------

\section{Fredholm theory of the Dirac operator} \label{fred}

We begin by recalling that points in the dual torus $\ksi\in\dual$
parametrises the set of flat holomorphic line bundles $L_\ksi\seta T$.
Moreover, such bundles have a natural choice of connection, denoted
$i\ksi$, which is consistent with the holomorphic structure; see
\cite{J2}.

In fact, $\dual$ also parametrises the set of flat holomorphic line
bundles over $\torus$. Using the projection $p_1:\torus\seta T$,
one obtains the holomorphic line bundle $p_1^*(L_{\ksi})$ over
$\torus$, which we shall also denote by $L_{\ksi}$ for simplicity;
let $\omega_{\ksi}$ be the pullback of the flat constant connection
on $L_{\ksi}\seta T$ described above; clearly, such connection is
also flat.

Let $E\seta\torus$ be a rank 2 bundle provided with an instanton
connection $A\in\cala_{(k,\ksi_0)}$. Form the bundle $E\otimes
L_{\ksi}$ with the corresponding connection $A_{\ksi}=A\otimes I +
I\otimes\omega_{\ksi}$; since all we have done was to add a flat
term to our original instanton, $A_{\ksi}$ is still an instanton on
the twisted bundle. We also require $A$ to be irreducible; clearly,
its twisted version $A_{\ksi}$ is also irreducible.

Consider now the Dirac operator acting on the bundle
$E(\ksi)=E\otimes L_{\ksi}$, coupled to the connection $A_{\ksi}$, and
its adjoint:
\begin{eqnarray*}
D_{A_{\ksi}} & : & \Gamma(E(\ksi)\otimes S^+)\seta\Gamma(E(\ksi)\otimes S^-) \\
D_{A_{\ksi}}^* & : & \Gamma(E(\ksi)\otimes S^-)\seta\Gamma(E(\ksi)\otimes S^+)
\end{eqnarray*}
where the spaces of sections are provided with norms suitably defined.
Since the base manifold is flat and the connection is anti-self-dual,
the Weitzenb\"ock formula on $E(\ksi)\otimes S^+\seta\torus$ is simply:
\begin{eqnarray}
D_{A_{\ksi}}^* D_{A_{\ksi}} & = & \nabla_{A_{\ksi}}^*\nabla_{A_{\ksi}}  \label{weit} \\
\Rightarrow \ \ ||D_{A_{\ksi}}s||^2 & = & ||\nabla_{A_{\ksi}}s||^2 \nonumber
\end{eqnarray}
Hence, if $A_{\ksi}$ is irreducible, there are no covariantly constant
sections of \linebreak $E(\ksi)\otimes S^+$. This means that the
kernel of  $D_{A_{\ksi}}$ is trivial. Now, if $D_{A_{\ksi}}$ is a
Fredholm operator, then  ${\rm ker}D_{A_{\ksi}}^*$ (which coincides
with ${\rm coker}D_{A_{\ksi}}$)  is a finite dimensional subspace of
$\Gamma(E(\ksi)\otimes S^-)$.

In this rather technical but fundamental section, we prove that this
is indeed the case:

\begin{thm} \label{dirac-fred}
Given any instanton connection $A\in\cfg$, the Dirac operators:
\begin{equation} \label{oopp}
D_{A_{\ksi}}^*:L^2_1(E(\ksi)\otimes S^-)\seta L^2(E(\ksi)\otimes S^+)
\end{equation}
form a smooth family of Fredholm operators parametrised by
$\dual\setminus\{\as\}$. Moreover, $\ind D_{A_\ksi}^*=k$, for all
$\ksi\in \dual\setminus\{\as\}$.
\end{thm}

The Sobolev norm in the left hand side of (\ref{oopp}) is defined as follows.
Let $D_\ksi^*$ be the Dirac operator $L_\ksi\otimes S^-\seta L_\ksi\otimes S^+$.
Then $L^2_1(E(\ksi)\otimes S^-)$ is the completion of $\Gamma(E(\ksi)\otimes S^-)$
with respect to the norm:
\begin{equation} \label{n.thm}
||s||_{L^2_1}=||s||_{L^2}+||D_\ksi^*s||_{L^2}
\end{equation}

The proof consists of three steps, which we now outline. We first
prove that the operators $D_\ksi^*:L^2_1(L_\ksi\otimes S^-)\seta L^2(L_\ksi\otimes S^+)$
are invertible for nontrivial $\ksi\in\dual$. A gluing argument then
shows that the Dirac operator coupled to a twisted instanton $A_\ksi$
is Fredholm if $\ksi\neq\ksi_0$. To compute the index, we use the
Gromov-Lawson Relative Index Theorem \cite{GL}.

It is important to note here that $D_{A_\ksi}$ fails to be
Fredholm when $\ksi=\as$; the reason will be clear from the
proof of the theorem. As we will see, this phenomenon is the source
of the singularities that appear in the transformed objects.

\subsection{The flat model.}

Let $L_\ksi\seta\torus$ be the flat line bundle described above, provided
with the constant connection $\omega_\ksi$. Our starting point to prove
theorem \ref{dirac-fred} is the following proposition.

\begin{prop} \label{flatmodel}
For non-trivial $\ksi\in\dual$, the coupled Dirac operator
$$ D_\ksi^*:L^2_1(L_\ksi\otimes S^-)\seta L^2(L_\ksi\otimes S^+) $$
is invertible. Its inverse is denoted by $Q^\infty_\ksi$.
\end{prop}

\pf
Let $L_{\ksi}\seta\torus$ be a flat line bundle as above, provided with
the constant connection $\omega_{\ksi}=p^*(-i\ksi)$, as described in
\cite{J2}. Consider the twisted Dirac operator:
$$ D_{\ksi}:\Gamma(L_{\ksi}\otimes S^+)\seta\Gamma(L_{\ksi}\otimes S^-) $$
and its adjoint $D_{\ksi}^*$.

Since $M=\torus$ is a K\"ahler surface, we have the
following decompositions:
\eq \label{decomp} \left\{ \begin{array}{l}
S^+=\Lambda^{(0,0)}_M L_{\ksi} \oplus \Lambda^{(0,2)}_M L_{\ksi} \\
S^-=\Lambda^{(0,1)}_M L_{\ksi} =\Lambda^{(0,1)}_{T} L_{\ksi} \oplus
    \Lambda^{(0,1)}_{\cpx}
\end{array} \right. \end{equation}
With respect to these decompositions, the Dirac operator and its adjoint are given by:
\eq \label{dirac} \begin{array}{cc}
D_{\ksi}=\matriz \del^{(z)}_{\ksi} & -\del^{(w),*}_{\ksi} \\
\del^{(w)}_{\ksi} & -\del^{(z),*}_{\ksi}  \matfim &
D^*_{\ksi}=\matriz -\del^{(z),*}_{\ksi} & -\del^{(w),*}_{\ksi} \\
\del^{(w)}_{\ksi} & \del^{(z)}_{\ksi}  \matfim
\end{array} \end{equation}
where $\del^{(z,w)}_{\ksi}$ denotes the Dolbeault operator twisted by $\omega_{\ksi}$ along the
toroidal ($z$) and plane ($w$) complex coordinates, i.e. the components of the covariant
derivative. Hence, the coupled Dirac laplacian $\triangle_{\ksi}=D_{\ksi}^*D_{\ksi}$ mapping
$\Lambda^{(0,0)}_M L_{\ksi}\oplus \Lambda^{(0,2)}_M L_{\ksi}$ to itself is just:
\eq \label{lap} \matriz
\triangle_{\ksi}^{(z)}+\triangle_{\ksi}^{(w)} & 0 \\
0 & \triangle_{\ksi}^{(z)}+\triangle_{\ksi}^{(w)}
\matfim \end{equation}
The off-diagonal terms are cancelled, for they are proportional to
the curvature, which was supposed to vanish. Moreover, the flat
connection $\omega_\ksi$ is a pull back from the torus, so that
$\triangle_{(w)_\ksi}$ is just the usual plane laplacian. Let us
concentrate on a single component, say $\Lambda^{(0,0)}_M L_{\ksi}$.

First, we want to solve the homogeneous equation $\triangle_{\ksi}f=0$ for \linebreak
$f\in\Lambda^{(0,0)}_M(L_{\ksi})$ and a fixed $\ksi\in\dual$.
Now, separate variables, supposing that $f(z,w)=\varphi(z)g(w)$:
$$ \triangle_{\ksi}f=0\ \see \ g\triangle_{\ksi}^{(z)}\varphi+\varphi\triangle^{(w)}g=0 $$
Therefore:
\eq \label{sep} \left\{ \begin{array}{l}
(i)\ \triangle_{\ksi}^{(z)}\varphi=\lambda^2\varphi \\
(ii)\ \triangle^{(w)}g=-\lambda^2 g\ \rightarrow\  (\triangle^{(w)}+\lambda^2)g=0
\end{array} \right. \end{equation}
where $\lambda^2$ are the eigenvalues of the $\ksi$-twisted laplacian over
the torus. They form a discrete, unbounded set $\{\lambda_n(\ksi)\}_{n\in\mathbb{N}}$
of $\real^+$, each being a function of the parameter $\ksi$. Note that
since $H^0(T,L_\ksi)=0$ for nontrivial $\ksi\in\dual$, we can indeed
guarantee that $\lambda_n(\ksi)>0$ for all nontrivial $\ksi$. On the
other hand, for $L_\ksi=\underline{\cpx}$, the laplacian has a
1-dimensional kernel, i.e. one zero eigenvalue.

As usual, we can decompose $f$ on the eigenstates of $\triangle_{\ksi}^{(z)}$,
i.e.:
\begin{equation} \label{f.decomp}
f=\sum_{n}g_n(w)\varphi_{n}(z)
\end{equation}
where $\{\varphi_n\}$ is an orthonormal basis for the $L^2$ norm
on $\Lambda^{(0,0)}_M(L_{\ksi})$ of eigenstates with eigenvalues
$\{\lambda^2_n\}$; so, $||f||_{L^2(\torus)}^2=\sum_{n}||g_n||^2_{L^2(\cpx)}$.
Moreover:
\eq \label{evalue}
\triangle_{\ksi}f=\sum_{n}[(\triangle^{(w)}+\lambda_n^2)g_n]\varphi_n
\end{equation}

\begin{prop} \label{est1}
Let $\rho\in L^2(L_{\ksi}\otimes S^+)$ be compactly supported and
suppose that $\ksi$ is nontrivial. Then there is $f\in L^2(L_{\ksi}\otimes S^+)$
and a constant $k<\infty$ such that $\Delta_{\ksi}f=\rho$ and
$||f||_{L^2}\leq k||\rho||_{L^2}$.
\end{prop}

\pf Given (\ref{evalue}), solving the equation $\triangle_{\ksi}f=\rho$
amounts to solve \linebreak $(\triangle^{(w)}+\lambda^2_n)g_n=\rho_n$
for each $n$, where $g_n,\rho_n$ are the components of $g,\rho$ along
the eigenspaces of $\lambda^2_n$, respectively.

Fix some integer $n$ and denote by $F_n$ the fundamental solution of
$(\Delta^{(w)}+\lambda^2_n)F_n(w)=0$. Rescale the plane coordinate
$w^{\prime}=\lambda_n w$, which transforms the previous equation to
$(\triangle^{(w^{\prime})}+1)F_n(\frac{w^\prime}{\lambda_n})=0$. The
unique integrable solution for this equation is the Bessel function
$K_0$ (see below), so that $F_n(w)=K_0(\lambda_n w)$. Solutions to
the non-homogeneous equations will then be given by the convolution:
\begin{equation} \label{conv}
g_n(w)=\int_{\real^{2}}F_n(w-x)\rho_n(x)dxd\overline{x}
\end{equation}
and recall that $||g_n||_{L^2}\leq||F_n||_{L^1}||\rho_n||_{L^2}$.
So, all we need is an estimate for $||F_n||_{L^1}$ which is independent
of $n$.

From the expression above, one sees that each $F_n$ is integrable if
the Bessel function $K_0$ is: $||F_n||_{L^1}=\lambda_n^{-2}||K_0||_{L^1}$.
So, let $\lambda={\rm min}\{\lambda_n\}_{n\in{\Bbb N}}$; therefore,
$||F_n||_{L^1}\leq\lambda^{-2}||K_0||_{L^1}$; putting
$k=\lambda^{-2}||K_0||_{L^1}$ we have $||g_n||_{L^2}\leq k||\rho_n||_{L^2}$
for each $n$. This completes the proof.  \pfend

Consider the Hilbert space $L^2_2(L_\ksi\otimes S^\pm)$ obtained by the completion of
$\Gamma(L_{\ksi}\otimes S^\pm)$ with respect to the norm:
\begin{equation} \label{norm}
||s||_{L^2_2}=||s||_{L^2}+||\triangle_{\ksi}s||_{L^2}
\end{equation}
The map $\triangle_{\ksi}:L^2_2(L_{\ksi}\otimes S^-)\seta L^2(L_{\ksi}\otimes S^-)$
is then bounded, for clearly $||\Delta_\ksi s||_{L^2}\leq||s||_{L^2_2}$.
Let $G_\ksi:L^2(L_{\ksi}\otimes S^-)\seta L^2_2(L_{\ksi}\otimes S^-)$
be the inverse of $\triangle_\ksi$ given by proposition
\ref{est1}. Using the inequality of the proposition, one shows that
$G_\ksi$ is also bounded, if $\ksi$ is nontrivial:
\begin{eqnarray*}
||G_\ksi s||_{L^2_2}& = & ||G_\ksi s||_{L^2}+||\triangle_\ksi G_\ksi s||_{L^2}\ = \
||G_\ksi s||_{L^2}+||s||_{L^2}\ \leq \\
& \leq & k ||s||_{L^2}+||s||_{L^2}\ \leq\ (k+1)\cdot||s||_{L^2}
\end{eqnarray*}
Moreover, we also conclude that:
\begin{equation} \label{gr.est}
||G_\ksi||< 1+\frac{C}{\lambda^2}
\end{equation}

Hence, $G_\ksi$ is an invertible operator when acting between the above
Hilbert spaces, if $\ksi$ is non-trivial.

\vskip10pt

\noindent {\bf Remarks:} We emphasise the necessity of assuming that
$\ksi$ is nontrivial. If \linebreak $\ksi=\hat{e}$, then the equation
(\ref{sep}i) admits one zero eigenvalue; on the other hand, the
fundamental solution of $\triangle^{(w)}g=0$ is essentially $\log r$,
which is not integrable. It is then impossible to get the estimate of
proposition \ref{est1}, in other words, the operator $\triangle_{(\ksi=\hat{e})}$
fails to be invertible. In addition, the parameter $k$ also depends on
$\ksi$, and $k\seta\infty$ (i.e. $\lambda\seta0$) as $\ksi\seta0$.

\vskip10pt

Now, define the norms:
\begin{equation} \label{metrics} \left\{ \begin{array}{l}
||s||_{L^2_1}=||s||_{L^2}+||D^*_{\ksi}s||_{L^2}\ {\rm if}\ s\in \Gamma(L_{\ksi}\otimes S^-)\\
||s||_{L^2_{l+1}}=||s||_{L^2_l}+||D_{\ksi}s||_{L^2_l}\ {\rm if}\ s\in \Gamma(L_{\ksi}\otimes S^+)
\end{array} \right. \end{equation}
and consider the Dirac operators as maps between the following Hilbert
spaces, obtained by the completion of $\Gamma(L_{\ksi}\otimes S^\pm)$
with respect to the above norms:
\begin{equation} \label{maps} \left\{ \begin{array}{l}
D^*_\ksi: L^2_1(L_{\ksi}\otimes S^-)\seta L^2(L_{\ksi}\otimes S^+)\\
D_\ksi: L^2_{l+1}(L_{\ksi}\otimes S^+)\seta L^2_l(L_{\ksi}\otimes S^-)
\end{array} \right. \end{equation}
Then $D^*_\ksi$ is clearly bounded. Furthermore, it has an inverse given by \linebreak
$(D^*_\ksi)^{-1}=D_\ksi G_\ksi:L^2(L_{\ksi}\otimes S^+)\seta L^2_1(L_{\ksi}\otimes S^-)$,
which is also bounded:
\begin{eqnarray*}
||(D^*_\ksi)^{-1}s||_{L^2_1} & = &
||(D^*_\ksi)^{-1}s||_{L^2}+||D^*_\ksi(D^*_\ksi)^{-1}s||_{L^2}=\ \\
& = &||D_\ksi G_\ksi s||_{L^2}+||s||_{L^2} \  = \
||D_\ksi G_\ksi s||_{L^2_1}\ \leq \\
& \leq & ||G_\ksi s||_{L^2_2}\ \leq\ (k+1)\cdot||s||_{L^2}
\end{eqnarray*}

So, $D^*_\ksi$ is also Fredholm when acting as in (\ref{maps}), and our proof is complete.
To further reference, we shall denote $Q_\ksi^\infty=(D^*_\ksi)^{-1}$;
note, moreover, that this is a bounded, elliptic, pseudo-differential operator of order $-1$.
\pfend

We are left with one point to establish: the integrability of the
fundamental solution of $(\triangle+1)F=0$ in the plane. Indeed, first note
that since the operator $\triangle+1$ has polar symmetry, then the
fundamental solution $F$ also has. After imposing this symmetry,
we obtain the following ODE, for $r>0$:
$$ (\triangle+1)F(r)=0 \imply F^{\prime\prime}+\frac{1}{r}F^\prime-F=0 $$
This is a Bessel equation with parameter $\nu=0$. Its
solutions are linear combinations of the Bessel functions of
imaginary argument $I_0$ and $K_0$ (see \cite{A}, chapter 11).
Below are possible integral representations for these functions:
\begin{eqnarray*}
K_0(r)=\int_{1}^{\infty}e^{-rt}(t^2-1)^{-\frac{1}{2}}dt
& \cite{GR} & 8.432.3 \\
I_0(r)=\int_{-1}^{1}\cosh(rt)(t^2-1)^{-\frac{1}{2}}dt
& \cite{GR} & 8.431.2
\end{eqnarray*}
It is easy to see that $I_0(r)$ increases exponentially with $r$; it is also
finite for $r=0$. For the purpose of finding a Green's function for the operator
$\triangle+1$, this solution can be eliminated.

With the help of a table of integrals, one finds out that $K_0$ is integrable; indeed:
$$ \int_{\real^2}K_0(r)d^2vol=\int_{0}^{\infty}\int^{2\pi}_{0}K_0(r)rdrd\theta=
2\pi\int_{0}^{\infty}rK_0(r)dr=2\pi $$
by \cite{GR} 6.561.16 (choosing $\mu=1$, $\nu=0$, $a=1$). This means that
\linebreak $||K_0||_{L^1}=2\pi$.

\begin{prop} \label{est2}
The solution $f$ of the flat laplacian problem $\Delta_{\ksi}f=\rho$
of proposition (\ref{est1}) decays exponentially if $\ksi$ is
nontrivial, in the sense that there is a real constant $\lambda>0$
such that:
$$ \lim_{r\seta\infty}e^{\lambda r}|f|<\infty $$
\end{prop}

\pf As $r\seta\infty$, the Bessel function $K_0$ admits the
following asymptotic expansion (\cite{Wa}, p.202):
\begin{equation} \label{k0exp}
K_0(r) \sim \left( \frac{\pi}{2} \right)^{\frac{1}{2}} \frac{e^{-r}}{\sqrt{r}}
              \left[ 1-\frac{1}{8r}+\frac{9}{128r^2}+\dots \right]
\end{equation}
Now since each $\rho_n$ has compact support, it follows from (\ref{conv})
that each $g_n$ will also decay exponentially:
$$ g_n(w) \sim \left( \frac{\pi}{2} \right)^{\frac{1}{2}}\cdot
\int_\Omega \frac{e^{-\lambda_n|w-x|}}{\sqrt{\lambda_n|w-x|}}
\left[ 1-\frac{1}{8\lambda_n|w-x|}+\dots \right]\rho_n(x)
dxd\overline{x} $$
where $\Omega$ is the support of $\rho$. As $|w|\seta\infty$, then also
$|w-x|\sim |w|$ for all $x\in\Omega$. Therefore,
$$ g_n(w) \sim \left( \frac{\pi}{2} \right)^{\frac{1}{2}}
\frac{e^{-\lambda_n|w|}}{\sqrt{\lambda_n|w|}}
\left[ 1-\frac{1}{8\lambda_n|w|}+\dots \right]
\cdot \int_\Omega \rho_n(x) dxd\overline{x},
\ \ {\rm as}\ \ |w|\seta\infty $$
Choosing $0<\lambda<{\rm min}\{\lambda_n\}_{n\in{\Bbb N}}$, the
statement follows from the eigenspace decomposition of $f$
(\ref{f.decomp}) and (\ref{evalue}).
\pfend

In particular, note that $(f/w)$ also belongs to $L^2(L_\ksi\otimes S^+)$.
Define \linebreak $\widetilde{L^2}(L_\ksi\otimes S^+)$ as the space of sections
$\psi$ such that $\psi/w$ is square-integrable. The proposition just proved
implies that the flat model laplacian acting as follows:
$$ \triangle_\ksi:\widetilde{L^2}(L_\ksi\otimes S^{\pm})\seta L^2(L_\ksi\otimes S^{\pm}) $$
is an invertible operator. Since $\triangle_\ksi=D_\ksi D^*_\ksi$, we conclude
that:
\begin{equation} \label{df2}
D_\ksi^*: \widetilde{L^2}(L_\ksi\otimes S^-) \seta L^2(L_\ksi\otimes S^+)
\end{equation}
is also invertible.

\subsection{Completing the proof of the theorem \ref{dirac-fred}.}
Let $K$ denote a closed ball in $\cpx$ of sufficiently large radius $R$;
its complement is $D_R$ defined as above. To show that $D_{A_{\ksi}}^*$
is Fredholm, first note that usual elliptic theory for compact manifolds
guarantees the existence of a parametrix for $D^*_{A_\ksi}$ inside this
compact core $T\times K$; this is a bounded, elliptic, pseudo-differential
operator:
$$ Q^K_{A_\ksi}:L^2(E(\ksi)\otimes S^+|_{T\times K})\seta
   L^2_1(E(\ksi)\otimes S^-|_{T\times K}) $$
of order $-1$.

On the other hand, it follows from lemma \ref{model} that:
$$ ||D^*_{A_\ksi}-
    (D^*_{\ksi_0+\ksi}\oplus D^*_{-\ksi_0+\ksi})||
    _{L^2(T\times D_R)}^2<2\epsilon $$
where $\epsilon$ can be made arbitrarily small. Thus,
$D^*_{A_\ksi}|_{T\times D_R}$ is also invertible for sufficiently
large $R\gg 0$, if $\ksi\neq\as$. Denote this inverse by
$Q^\infty_{A_\ksi}$; this is also a bounded, elliptic, pseudo-differential
operator of order $-1$.

Now choose $\beta_1,\beta_2:\cpx\seta\real$ respectively supported over
$K$ and $D_R$ and satisfying $\beta_1^2+\beta_2^2=1$ everywhere.
We can patch together our two parametrix $Q^K_{A_\ksi}$ and $Q^{\infty}_{A_\ksi}$
in the following way:
\eq \label{parametrix}
P_{A_\ksi}g=\beta_1Q^K_{A_\ksi}(\beta_1g)+\beta_2Q^\infty_{A_\ksi}(\beta_2g)
\end{equation}
This is the same as restricting the section $g$ to $T\times K$
(respectively, \linebreak $T\times D_R$), apply
$Q^K_{A_\ksi}$ ($Q^\infty_{A_\ksi}$) and restricting the result again
to $T\times K$ ($T\times D_R$). Note that $P_{A_\ksi}$ acts as follows:
$$ P_{A_\ksi}:L^2(E(\ksi)\otimes S^+)\seta L^2_1(E(\ksi)\otimes S^-). $$

We want to show that this is a parametrix for $D_{A_{\ksi}}^*$.
In fact, take \linebreak $g\in L^2(E(\ksi)\otimes S^+)$; then:
\begin{eqnarray}
D_{A_{\ksi}}^*P_{A_\ksi}g &=& D_{A_{\ksi}}^*[\beta_1Q^K_{A_\ksi}(\beta_1g)]+
D_{A_{\ksi}}^*[\beta_2Q^\infty_{A_\ksi}(\beta_2g)]= \nonumber \\
& = & \{\beta_1D_{A_{\ksi}}^*Q^K_{A_\ksi}(\beta_1g)+
\beta_2D_{A_{\ksi}}^*Q^\infty_{A_\ksi}(\beta_2g)\}+ \label{xx} \\
& & +\underbrace{d\beta_1.Q^K_{A_\ksi}(\beta_1g)+
d\beta_2.Q^\infty_{A_\ksi}(\beta_2g)} \nonumber \\
& & \hspace{6.5em} S^\infty g \nonumber
\end{eqnarray}
where ``$.$'' means Clifford multiplication.

Since $Q^K_{A_\ksi}$ is a parametrix for $D_{A_{\ksi}}^*$ inside $T\times K$,
the first term (inside brackets) equals the identity plus a compact
operator $S^K$ acting on $\beta_1g$. Similarly, in the second term,
$Q^{\infty}_{A_\ksi}$ is the inverse of the Dirac operator outside $K$.
Together, the first two terms form the identity operator plus $S^K$. Hence:
$$ (D_{A_\ksi}^*P_{A_\ksi}-I)g=S^Kg+S^\infty g $$
where $S^\infty:L^2(E(\ksi)\otimes S^+)\seta L^2(E(\ksi)\otimes S^+)$
is the operator over the brackets in
(\ref{xx}). Since $Q^K_{A_\ksi}$ and $Q^\infty_{A_\ksi}$ are bounded operators,
so is $S^\infty$; we argue that this is also a compact operator.

In fact, let $\widetilde{\partial K}$ denote the closure of the the support of
$d\beta_1$ and $d\beta_2$ (which is an annulus around the boundary of $K$).
Consider the diagram:
$$ \begin{array}{ccc}
L^2(E(\ksi)\otimes S^+)&\stackrel{s}{\longrightarrow}&
L^2_1(E(\ksi)\otimes S^+|_{T^2\times\widetilde{\partial K}})\\
& &\downarrow\ i \\
& &L^2(E(\ksi)\otimes S^+|_{T^2\times\widetilde{\partial K}})\\
& & \cap \\ & &L^2(E(\ksi)\otimes S^+)
\end{array} $$
Now, let $\Upsilon\subset L^2(E(\ksi)\otimes S^+)$ be a
bounded set; since $s$ is a bounded operator, $s(\Upsilon)$ is also
bounded. By the Rellich lemma, the map $i$ is a compact inclusion;
note that $\widetilde{\partial K}$ is a compact subset of the plane.
Hence, $i(s(\Upsilon))$ is a relatively compact subset of
$L^2(E(\ksi)\otimes S^+|_{T^2\times\widetilde{\partial K}})$, and
clearly also a relatively compact subset of $L^2(E(\ksi)\otimes S^+)$.
This means that:
$$ S^\infty=i\circ s:L^2(E(\ksi)\otimes S^+)\seta L^2(E(\ksi)\otimes S^+) $$
is a compact operator, as have we claimed. We conclude that:
$$ D_{A_{\ksi}}^*P_{A_\ksi}-I=[{\rm compact\ operator}] $$
so that (\ref{parametrix}) is indeed a parametrix for $D_{A_\ksi}^*$ if
$\ksi\neq\pm\ksi_0$.

Finally, to compute the index of $D_{A_\ksi}^*$ we use the
Relative Index Theorem of Gromov \& Lawson \cite{GL} 
(see also the appendix in \cite{J}). One can show that:

\begin{lem} \label{indmod}
If $A\in\cala_{(k,\ksi_0)}$, then ${\rm index}D^*_{A_\ksi}=k$.
\end{lem}

\paragraph{The Green's operator.}
Clearly, the Dirac laplacian, with the norms as in (\ref{norm}):
\begin{eqnarray}
\Delta_{A_\ksi}: & L^2_2(E\otimes L_\ksi\otimes S^+)\seta
                  L^2(E\otimes L_\ksi\otimes S^+) & \nonumber \\
& \Delta_{A_\ksi}=D_{A_\ksi}^*D_{A_\ksi} & \label{green}
\end{eqnarray}
is also a Fredholm operator. In particular, by general Fredholm
theory, there is a bounded operator $G_{A_\ksi}$, called the Green's
operator, such that:
$$ \Delta_{A_\ksi}G_{A_\ksi}=Id-H_\ksi $$
where $H_\ksi$ is the finite rank orthogonal projection operator:
$$ H_\ksi:L^2_2(E\otimes L_\ksi\otimes S^+)\seta{\rm ker}(\Delta_{A_\ksi}) $$

\subsection{Harmonic spinors and cohomology.} \label{spch}
To conclude this chapter, we want to interpret the harmonic spinors \linebreak
$\psi\in{\rm ker}D_A^*$ as some holomorphic object defined in terms of the
compactified bundle $\ee\seta\tproj$. Indeed, we aim to establish the
following identification:

\begin{prop} \label{spin/coho}
If $A$ has nontrivial asymptotic state $\ksi_0\in\dual$ and $k>0$,
then there is an isomorphism $H^1(\tproj,\ee)\equiv{\rm ker}D_A^*$.

\end{prop}

Note that ${\rm ker}D_A^*\subset L^2_1(E\otimes S^-)$, with the norm
defined in (\ref{n.thm}). First, we must show that $H^1(\tproj,\ee)$ has
the correct dimension.

\paragraph{Vanishing theorem.}
Since $\chi(\ee)=-k$, it is enough to show that the cohomologies of
orders $0$ and $2$ vanish in order to conclude that \linebreak $h^1(\tproj,\oo(\ee))=k$.

Let us assume that the restriction of $\ee$ to the elliptic
curves $\ee|_{T\times\{w\}}$ is semi-stable for all $w\in\proj$.
We can then regard $\ee\seta\tproj$ as a family of extensions:
$$ 0 \seta L_\ksi \seta \ee|_{T_w} \seta L_{-\ksi} \seta 0 $$
of a flat line bundle $L_\ksi$ by its dual $L_{-\ksi}$, where
$\ksi\in\dual$ depends holomorphically on $w\in\proj$; in other
words, the family is parametrised by $\proj$.

Since such extension bundles can be indecomposable if and only if
$\ksi=-\ksi$ (i.e. $\ksi$ has order 2 in $\dual$), we conclude that
$\ee|_{T_w}$ splits as a sum of flat line bundles for all but
finitely many points $w\in\proj$. Furthermore, these flat line
bundles are holomorphically nontrivial for all but finitely many
points $w\in\proj$.

This observation leads to the desired vanishing result:

\begin{lem} \label{vanish}
If $\ee$ is irreducible and $k>0$, then:
$$ h^0(\tproj,\ee(\xi))=h^2(\tproj,\ee(\xi))=0, \ \forall\ksi\in\dual $$
\end{lem}

Let $L_\ksi\seta T$ be a flat line bundle as described in \cite{J2};
denote:
$$ \ee(\ksi)=\ee\otimes p_1^*L_\ksi \ \ \ {\rm and} \ \ \
   \tilde{\ee}(\ksi)=\ee\otimes p_1^*L_\ksi\otimes p^*_2\oo_{\proj}(1) $$
Note that we can regard $p^*_2\oo_{\proj}(1)$ as the line bundle corresponding
to the divisor $T_\infty$. It follows from the lemma that:
$$ h^1(\tproj,\ee(\ksi))=h^1(\tproj,\tilde{\ee}(\ksi))=k $$
for every $\ksi\in\dual$.

\pf Take $w\in\proj$ such that $\ee(\ksi)|_{T_w}=L_{\ksi_1}\oplus L_{\ksi_2}$
for some non-trivial $\ksi_1,\ksi_2\in\dual$; the existence of such
point follows from the observations made prior to the statement of
the lemma. Let $V\subset\proj$ be an open neighbourhood of $w$ such
that every point of $V$ satisfy a the same condition.

Suppose there is a holomorphic section $s\in H^0(M,\ee(\ksi))$; it
gives rise to a holomorphic section $s_w$ of $\ee(\ksi)|_{T_w}\seta
T_w$. On the other hand, we have that \linebreak
$h^0(T,\ee(\ksi)|_{T\times\{w\}})=0$, hence $s_w\equiv0$. Moreover,
$s_w\equiv0$ for all $w\in V$, so that $s$ must vanish identically on
the open set $T\times V$, hence vanish everywhere and
$h^0(\ee(\ksi))=0$. The vanishing of $h^0(\tilde{\ee}(\ksi))$ is proved
in the very same way by noting $\tilde{\ee}(\ksi)|_{T_w}\equiv \ee(\ksi)|_{T_w}$
since $p^*_2\oo_{\proj}(1)|_{T_w}=\underline{\cpx}$.

The vanishing of the $h^2$'s follows from Serre duality and a similar argument for the
bundle $\ee(\ksi)\otimes K_{\proj}$. More precisely, Serre duality
implies that:
$$ \begin{array}{rcl}
H^2(\tproj,\ee(\ksi)) & = & H^0(\tproj,\ee(\ksi)^\vee\otimes K_{\tproj})^* \\
& = & H^0(\tproj,\ee(\ksi)^\vee\otimes p_2^*\oo_{\proj}(-2))^*
\end{array} $$
On the other hand, it is easy to see that:
$$ \ee(-\ksi)|_{T_w} \equiv (\ee(\ksi)^\vee\otimes p_2^*\oo_{\proj}(-2))|_{T_w} $$
so that we can apply the same argument as above to show that
\linebreak $h^0(\tproj,\ee(\ksi)^\vee\otimes K_{\tproj})=0$.
\pfend

\paragraph{Proof of proposition \ref{spin/coho}.}
Let $\{w_i\}$ be the set of points in $\proj$ for which
$H^0(T_{w_i},\ee|_{T_{w_i}})\neq\{0\}$. As we argued above,
there are only finitely many such points; in fact, it can be
shown that there are at most $k$ such points (see lemma 2 of
\cite{J2}). Suppose that $\#\{w_i\}=p\leq k$; note also that
$\infty\notin\{w_i\}$ if $\ksi_0$ is nontrivial.

Denote by $B$ the divisor in $\tproj$ consisting of the elliptic
curves lying over these points, i.e. $B=\sum_i T\times\{w_i\}$. Also,
denote $\ee(p)=\ee\otimes\oo_{\tproj}(B)$.

Consider the exact sequence of sheaves:
$$ 0\seta\oo(\ee)\seta\oo(\ee(p))\seta\oo(\ee(p)|_B)\seta 0 $$
which induces the following sequence of cohomology:
\footnotesize
\begin{equation} \label{xy}  \begin{array}{ccccc}
0 \seta H^0(B,\ee(p)|_B) \seta & \underbrace{H^1(\tproj,\ee)}
  & \seta & \underbrace{H^1(\tproj,\ee(p))} & \seta H^1(B,\ee(p)|_B) \seta 0 \\
 & {\rm dim}=k & & {\rm dim}=k &
\end{array} \end{equation}
\normalsize \baselineskip18pt
and note that $p\leq h^0(B,\ee(p)|_B)=h^1(B,\ee(p)|_B)\leq2k$.
It follows from (\ref{xy}) that $h^0(B,\ee(p)|_B)=h^1(B,\ee(p)|_B)=k$,
so that the map $H^0(B,\ee(p)|_B)\seta H^1(\tproj,\ee)$ is an
isomorphism.

This means that each element in $H^1(\tproj,\ee)$ can be represented
by a $(0,1)$-form $\theta$ supported on tubular neighbourhoods of the
fibres $T\times\{w_i\}$. Pulling $\theta$ back to $\torus$, we obtain
a compactly supported $(0,1)$-form, which we also denote by $\theta$,
since $\ksi_0$ is nontrivial.

We want to fashion a solution $\psi$ of $D_A^*\psi=0$ out of
$\theta$, and within the same cohomology class. In other words,
we want to find a section $s\in L^2(\Lambda^0E)$ such that
$D_A^*(\theta+\del_As)=0$. Since $D_A^*=\del_A^*-\del_A$, this
is the same as solving the equation:
$$ \del_A^*\del_As=\Delta_As=-\del_A^*\theta $$
for a compactly supported $\theta$.

In the Fredholm theory for the Dirac operator developed above, we
constructed the Green's operator $G_A$ of the Dirac laplacian
$\Delta_A$. Thus, we can write $s=-G_A\del_A^*\theta$ and
$\psi=\theta-\del_AG_A\del_A^*\theta=P\theta$, where $P$ denotes
the $L^2$ projection $L^2(E\otimes S^-)\stackrel{P}{\seta}{\rm ker}D^*_A$.

We must verify that $\psi\in L^2(E\otimes S^-)$; it is enough to
show that $\del_AG_A\del_A^*\theta$ is square-integrable for any
compactly supported (0,1)-form $\theta$. First note that
$\gamma=\del_A^*\theta$ also has compact support, thus
$s=G_A\gamma\in L^2(\Lambda^0E)$. Therefore, we have:
\begin{eqnarray*}
||\del_A s||_{L^2}^2 & = & <\del_A s,\del_A s> \ = \
<\del_A s,(\del_AG_A)\gamma> \ = \\
& = & <(\del_AG_A)^*\del_A s,\gamma>
\end{eqnarray*}
which is finite, since $\gamma$ is compactly supported. Note the the integration
by parts made from the first to the second line is justified by the same fact.
Therefore, $\psi$ is indeed a square-integrable solution of $D_A^*\psi=0$.

Finally, to see that the map defined above is injective (hence an isomorphism),
let $\theta'$ be another $(0,1)$-form supported around $B$ and within
the same cohomology class as $\theta$, so that $\theta-\theta'=\del_A\alpha$.
Thus:
\begin{equation} \begin{array}{rcl}
\psi-\psi' & = & (\theta-\del_AG_A\del_A^*\theta)-(\theta'-\del_AG_A\del_A^*\theta') = \\
& = & (\theta-\theta')-\del_AG_A\del_A^*(\theta-\theta') = \\
& = & \del_A\alpha-\del_AG_A\del_A^*\del_A\alpha = \del_A\alpha-\del_A\alpha \ = \ 0
\end{array} \end{equation}

This completes the proof. \pfend

%--------------------------------------------------------------------

\section{Nahm transform of doubly-periodic \newline instantons} \label{nahm}

Recall that our starting point is a rank two vector bundle \linebreak
$E\seta\torus$ provided with an instanton connection
$A\in\cala_{(k,\ksi_0)}$, where the instanton number $k$ and the
asymptotic state $\ksi_0$ are from now on fixed.

Over the punctured Jacobian torus $\dual\setminus\{\as\}$, consider
the trivial Hilbert bundle $\hat{H}\seta\dual\setminus\{\as\}$ whose fibres are
$\hat{H}_{\ksi}=L^2_1(E(\ksi)\otimes S^-)$. Taking the $L^2_1$-norm
on the fibres, $\hat{H}$ becomes an hermitian bundle. Moreover,
call $\hat{d}$ the trivial connection on $\hat{H}$; such connection
is clearly unitary, hence one can define a holomorphic structure
over $\hat{H}$.

Now consider the finite-dimensional sub-bundle $V\hookrightarrow\hat{H}$
over $\dual\setminus\{\as\}$ whose fibres are given by
$V_{\ksi}={\rm ker}D_{A_{\ksi}}^*$. Remark that this is actually the
{\em index bundle} for the family of Dirac operators $D_{A_{\ksi}}$.
Let $i:V\seta\hat{H}$ be the natural inclusion and $P:\hat{H}\seta V$
the fibrewise orthogonal $L^2$ projection; more precisely,
$P_\ksi=I-D_{A_\ksi}G_{A_\ksi}D^*_{A_\ksi}$ for each
$\ksi\in\dual\setminus\{\as\}$, where $G_{A_\ksi}$ denotes the
Green's operator for (\ref{green}), $I$ is the identity operator.
We can define a connection on $V$ via the projection formula:
\eq \label{ft}
\nabla_{B}=P\circ\hat{d}\circ i
\end{equation}
where $B$ is the associated connection form.

Clearly, $V$ inherits the hermitian metric from $\hat{H}$, and $B$
is also unitary with respect to this induced metric. Hence, we can
provide $V$ with the holomorphic structure coming from the unitary
connection $B$.

Alternatively, $V$ also admits an interpretation in terms of {\em monads},
see \cite{DK}. The Dirac operator can be unfolded into a family of elliptic
complexes parametrised by $\dual\setminus\{\as\}$, namely:
\begin{equation} \label{monad}
0 \seta L^2_2(\Lambda^0E(\ksi)) \stackrel{\del_{A_\ksi}}{\longrightarrow}
   L^2_1(\Lambda^{0,1}E(\ksi)) \stackrel{-\del_{A_\ksi}}{\longrightarrow}
   L^2(\Lambda^{0,2}E(\ksi)) \seta 0
\end{equation}
which, of course, are also Fredholm. Moreover, the cohomologies of
order 0 and 2 must vanish, by proposition \ref{vanish}. As in \cite{DK},
such holomorphicfamily defines a holomorphic vector bundle
$V\seta(\dual\setminus\{\as\})$, with fibres
$V_\ksi=H^1(\ksi)={\rm ker}D_{A_\ksi}^*$, plus an unitary connection,
induced by orthogonal projection, which is compatible with the given
holomorphic structure. Such connection will be denoted by $B$. We will
invoke this construction repeatedly throughout this work.

The curvature $F_B$ of $B$ is simply given by:
$$ F_B=\nabla_B\nabla_B=P\hat{d}(P\hat{d}) $$
Explicit formulas for the matrix elements on an arbitrary local
trivialisation of $V\seta(\dual\setminus\{\as\})$ will be useful
later on. For instance, pick up an orthonormal frame
$\{\psi_i\}_{n=1}^{k}$ over an open set $U\subset\dual\setminus\{\as\}$.
Then, we have that:
\begin{eqnarray}
(B)_{ij} & = & <\psi_j,\nabla_B\psi_i> \ = \ <\psi_j,\hat{d}\psi_j> \nonumber \\
(F_B)_{ij} & = & <\psi_j,F_B\psi_i> \ = \ <\psi_j,P\hat{d}(P\hat{d}\psi_i)> \ =
\ <\psi_j,\hat{d}(P\hat{d}\psi_i)> \label{concurv}
\end{eqnarray}

\paragraph{Higgs field.}
We now define the Higgs field $\Phi\in{\rm End}(V)\otimes K_{\dual}$.
Let $w$ is the complex coordinate of the plane, and $\psi\in\Gamma(V)$,
i.e. for each $\ksi\in\dual\setminus\{\as\}$, $\psi[\ksi]\in{\rm ker}D^*_{A\ksi}$.
For a fixed $\ksi'$, the Higgs field will act on $\psi[\ksi']$ by
multiplying this section by the plane coordinate $w$ and then projecting i
t back to ${\rm ker}D^*_{A_\ksi}$:
\begin{equation} \label{higgs.gt}
(\Phi(\psi))[\ksi'] \ = \ 2\sqrt{2}\pi P_{\ksi'}(w\psi[\ksi'])d\ksi
\end{equation}
Its conjugate is clearly given by
$(\Phi^*(\psi))[\ksi'] = 2\sqrt{2}\pi  P_{\ksi'}
(\overline{w}\psi[\ksi'])d\overline{\ksi}$

Again, there is a subtle analytical point here. The spinors $\psi$
belong to $L^2(E(\ksi)\otimes S^-)$ but is not necessarily the case
that $w\psi$ also belong to $L^2(E(\ksi)\otimes S^-)$. However, we
have the following lemma:
\begin{lem}
If $\psi\in{\rm ker}D_A^*$ and $A$ has nontrivial asymptotic state, then
\linebreak $w\psi\in L^2(E\otimes S^-)$.
\end{lem}
\pf
The key result here is proposition \ref{est2}, and the observation
that follows it, in particular the invertibility of the operator
(\ref{df2}).

Let $K\subset\torus$ be a compact subset such that $D_A^*$ is
sufficiently close to the flat Dirac operator $D_{\as}^*$ outside
$K$. Thus, restricted to the complement of $K$, $D_A^*$ is
invertible acting from $\tilde{L^2}\seta L^2$.

Now if $\psi\in{\rm ker}D_A^*$, then $D_A^*(w\psi)=dw\cdot\psi\in
L^2(E(\ksi)\otimes S^+|_{\torus\setminus K})$ and the proposition
follows.
\pfend

Note that the dependence of $(B,\Phi)$ on the original instanton
$A$ is contained on the $L^2$-projection operator $P$, i.e. on the
$k$ solutions of $D_{A_\ksi}^*\psi=0$. It is easy to see that the
finite dimensional space spanned by these $\psi$ is gauge
invariant; moreover, the multiplication by $w$ also commutes with
gauge transformations $\hat{g}\in{\rm Aut}(V)$. Therefore, we have
that:

\begin{prop} \label{eqv2}
If $A$ and $A'$ are gauge equivalent irreducible instantons, then
the corresponding pairs $(B,\Phi)$ and $(B^\prime,\Phi^\prime)$
are also gauge equivalent.
\end{prop}

A pair $(B,\Phi)$ is called a {\em Higgs pair} on the bundle
$V\seta\dual\setminus\{\as\}$ if it satisfies Hitchin's self-duality
equations:
\eq \label{hiteq2}
\left\{ \begin{array}{l}
{\rm (i)}\ F_{B}+[\Phi,\Phi^*]=0 \\
{\rm (ii)}\ \overline{\partial}_{B}\Phi=0
\end{array} \right. \end{equation}

Recall from \cite{J2} that the unitary connection of the Poincar\'e
line bundle ${\bf P}\seta T\times\dual$ and its corresponding curvature
are given by:
$$ \omega(z,\ksi)=i\pi\cdot\sum_{\mu=1}^{2} \big(\ksi_\mu dz_\mu - z_\mu d\xi_\mu\big)
   \ \ {\rm and} \ \
   \Omega(z,\ksi)=2i\pi\cdot\sum_{\mu=1}^{2}d\ksi_\mu\wedge dz_\mu $$
From Braam \& Baal \cite{BVB}, we know that if $s\in\Gamma(E(\ksi)\otimes S^-)$,
then:
$$ D^*_{A_\ksi}(\hat{d}s)=[D^*_{A_\ksi},\hat{d}]s=-\Omega\cdot s $$
where ``$\cdot$" means Clifford multiplication. The local formula for the
curvature (\ref{concurv}) may now be cast on a more convenient form:
\begin{eqnarray*}
(F_B)_{ij} & = & <\psi_j,\hat{d}(P\hat{d}\psi_i)> \ = \
<\psi_j,\hat{d}(D_{A_\ksi}G_{A_\ksi}D^*_{A_\ksi}\hat{d}\psi_i)> \ = \\
& = & <-D^*_{A_\ksi}\hat{d}\psi_j,G_{A_\ksi}(D^*_{A_\ksi}\hat{d}\psi_i)> \ = \
<\Omega.\psi_j,G_{A_\ksi}(\Omega\cdot\psi_i)>
\end{eqnarray*}
Since the Clifford multiplication commutes with the Green's operator, we end up with:
\eq \label{curv} \begin{array}{rcl}
(F_B)_{ij}& = & -<(\Omega\wedge\Omega)\cdot\psi_i,G_{A_\ksi}\psi_i> \ = \\
& = & 8\pi^2 <(dz_1\wedge dz_2)\cdot\psi_j,G_{A_\ksi}\psi_i>d\ksi_1\wedge d\ksi_2 \ = \\
& = & -4\pi^2i<(dz_1\wedge dz_2)\cdot\psi_j,G_{A_\ksi}\psi_i>d\ksi\wedge d\overline{\ksi}
\end{array} \end{equation}
Note moreover that the inner product is taken in $L^2(E(\ksi)\otimes S^-)$,
integrating out the $(z,w)$ coordinates.

\begin{thm} \label{soln}
If $A\in\cfg$, then the associated pair $(B,\Phi)$ on the dual bundle
$V\seta\dual\setminus\{\as\}$ constructed above satisfies the Hitchin's
equations (\ref{hiteq2}).
\end{thm}

\pf Choose an open set $U\in\dual\setminus\{\as\}$ and pick up a
local orthonormal trivialisation of $V\seta\dual\setminus\{\as\}$
over $U$, such that the corresponding local frame $\{\psi_i\}_{n=1}^{k}$
is parallel at $\ksi$. Recall that $\psi_i(\ksi)\in{\rm ker}D^*_{A_\ksi}$.

First, we shall look at the second equation of (\ref{hiteq2}), and recall that \linebreak
$\dual\setminus\{\as\}$ was given the flat Euclidean metric induced from the quotient.
Once a local trivialisation is chosen, the endomorphism $\Phi$ can then be put
in matrix form, with matrix elements given by:
$$ a_{ij}(\ksi)=<\psi_j(\ksi),\Phi[\psi_i](\ksi)> $$
where $<,>$ is the inner product on $L^2(E(\ksi)\otimes S^-)$,
integrating out the $(z,w)$ coordinates. Clearly, $\Phi$ is a
holomorphic endomorphism if its matrix elements in a holomorphic
trivialisation are holomorphic functions. However:
$$ \Phi[\psi_i](\ksi)=P_\ksi(w\psi_i(\ksi))d\overline{\ksi}=
(I-D_{A_\ksi}G_{A_\ksi}D^*_{A_\ksi})(w\psi_i(\ksi))d\overline{\ksi} $$
so that:
\begin{eqnarray*}
a_{ij}(\ksi) & = & \! 2\sqrt{2}\pi \left\{<\psi_j(\ksi),w\psi_i(\ksi)>-
<\psi_j(\ksi),D_{A_\ksi}G_{A_\ksi}D^*_{A_\ksi}(w\psi_i(\ksi))>\right\}= \\
& = & \! 2\sqrt{2}\pi  \left\{<\psi_j(\ksi),w\psi_i(\ksi)>-
<D^*_{A_\ksi}\psi_j(\ksi),G_{A_\ksi}D^*_{A_\ksi}(w\psi_i(\ksi))>\right\}= \\
& = & \! 2\sqrt{2}\pi <\psi_j(\ksi),w\psi_i(\ksi)>
\end{eqnarray*}
Therefore:
\begin{eqnarray*}
\frac{\partial a_{ij}}{\partial\overline{\ksi}}(\ksi) &=&
2\sqrt{2}\pi  \left\{ \langle \partial_{B}\psi_j,w\psi_i\rangle +
\langle\psi_j,\del_{B}(w\psi_i) \rangle \right\} \ =
\\ & = & 2\sqrt{2}\pi \langle \psi_j,
         \left(\frac{\partial w}{\partial\overline{\ksi}}\right)\psi_i+
         \del_{B}\psi_i \rangle \ = \ 0
\end{eqnarray*}
as $\psi_i$ is parallel at $\ksi$. Since this can be done for all
$\ksi\in\dual\setminus\{\as\}$, the second equation is satisfied.

Now, we move back to (\ref{hiteq2}(i)). Let us first
compute the matrix elements $([\Phi,\Phi^*])_{ij}$. Note that:
\eq \label{ids} \left\{ \begin{array}{l}
(i)\ [D_{A_\ksi}^*,\overline{w}]\psi_i(\ksi)=D_{A_\ksi}^*(\overline{w}\psi_i(\ksi))=
-d\overline{w}\cdot\psi_i(\ksi)\\
(ii)\ [D_{A_\ksi}^*,w]\psi_i(\ksi)=D_{A_\ksi}^*(w\psi_i(\ksi))=0
\end{array} \right. \end{equation}
where we used the fact that $D_{A_\ksi}=\del^*_{A_\ksi}-\del_{A_\ksi}$.

Recall that for 1-forms $[\Phi,\Phi^*]=\Phi\Phi^*+\Phi^*\Phi$.
We compute each term separately:
\begin{eqnarray*}
\Phi^*\Phi(\psi_i) & = &  8\pi^2 P[\overline{w}P(w\psi_i)]d\ksi\wedge d\overline{\ksi}=\\
&=&  8\pi^2  \left\{\overline{w}P(w\psi_i)-D_{A_\ksi}G_{A_\ksi}
D^*_{A_\ksi}\overline{w}P(w\psi_i)\right\}d\ksi\wedge d\overline{\ksi}=\\
&=&  8\pi^2  \left\{\overline{w}w\psi_i-
\overline{w}D_{A_\ksi}G_{A_\ksi}D^*_{A_\ksi}(w\psi_i)-\right.\\
& &\left.-D_{A_\ksi}G_{A_\ksi}D^*_{A_\ksi}\overline{w}P(w\psi_i)\right\}
d\ksi\wedge d\overline{\ksi} \\
\Phi\Phi^*(\psi_i) & = &  8\pi^2  P[wP(\overline{w}\psi_i)]d\overline{\ksi}\wedge d\ksi=\\
&=&  8\pi^2  \left\{w\overline{w}\psi_i-wD_{A_\ksi}G_{A_\ksi}
D^*_{A_\ksi}(\overline{w}\psi_i)-\right.\\
& &\left.-D_{A_\ksi}G_{A_\ksi}D^*_{A_\ksi}wP(\overline{w}\psi_i)\right\}
d\overline{\ksi}\wedge d\ksi
\end{eqnarray*}

The two first terms of $\Phi\Phi^*$ and $\Phi^*\Phi$ cancel each other
and the third terms will cancel out when we take the inner
product with $\psi_j$. Moreover, the second term of $\Phi^*\Phi$ is
zero by (\ref{ids}(ii)). So we are left with:
\begin{eqnarray*}
([\Phi,\Phi^*])_{ij} & = & 8\pi^2 <\psi_j,[\Phi,\Phi^*]\psi_i> = \\
& = & 8\pi^2 <\psi_j,wD_{A_\ksi}G_\ksi
D^*_{A_\ksi}(\overline{w}\psi_i)>d\ksi\wedge d\overline{\ksi} = \\
&=&  8\pi^2  <D^*_{A_\ksi}(\overline{w}\psi_j),G_\ksi
D^*_{A_\ksi}(\overline{w}\psi_i)> d\ksi\wedge d\overline{\ksi} = \\
&=&-8\pi^2  <(dw\wedge d\overline{w}).\psi_j,G_\ksi\psi_i>d\ksi\wedge d\overline{\ksi} = \\
&=& - 4\pi^2 i<(dw_1\wedge dw_2).\psi_j,G_\ksi\psi_i>d\ksi\wedge d\overline{\ksi}
\end{eqnarray*}
where we have once more used the fact that the Clifford multiplication commutes with the
Green's operator. Summing the final expression above with ({\ref{curv}), one gets:
$$ (F_B)_{ij}+([\Phi,\Phi^*])_{ij} = - 4\pi^2 i <(dz_1\wedge dz_2+dw_1\wedge dw_2)\cdot\psi_j,G_\ksi\psi_i>
d\ksi\wedge d\overline{\ksi}=0 $$
for the first term of the inner product is zero since it consists of a
self-dual form (the K\"ahler form) acting on a negative spinor.
\pfend

Clearly, the above result has two weak points: it tells nothing
about the behaviour of the Higgs field around the singular points
$\as$; and it fails to show that the Higgs pairs so obtained are
admissible in the sense of \cite{J2}. In fact, establishing the
first point requires the use of algebraic-geometric methods, and
will be taken up in section \ref{holo} below. The second point will
be clarified in section \ref{inv}.

%--------------------------------------------------------------------

\section{Holomorphic version} \label{holo}

The vanishing results of section \ref{spch} put us in position to
define the transformed bundle $\vv\seta\dual$. Indeed, consider the
following elliptic complex:
\begin{equation} \label{cpx1}
0 \seta L^2_2(\Lambda^0\ee(\ksi)) \stackrel{\del_{A_\ksi}}{\seta}
L^2_1(\Lambda^{0,1}\ee(\ksi)) \stackrel{-\del_{A_\ksi}}{\seta}
L^2(\Lambda^{0,2}\ee(\ksi)) \seta 0
\end{equation}
According to proposition \ref{vanish}, $H^1(\tproj,\ee(\ksi))$ is
the only nontrivial cohomology of this complex. It then follows
that the family of vector spaces given by
$\vv_\ksi=H^1(\tproj,\ee(\ksi))$ forms a holomorphic vector bundle
of rank $k$ over $\dual$; denote such holomorphic structure by
$\del_{\vv}$. Note that $\vv_\ksi$ is defined even if $\ksi=\as$.
Furthermore, by proposition \ref{spin/coho}, $\vv|_{\dual\setminus\as}$
coincides holomorphically with the dual bundle $V$ defined on the
previous section, i.e.:
$$ (\vv,\del_{\vv})|_{\dual\setminus\{\as\}}\simeq(V,\del_B) $$

Moreover, $\vv$ comes equipped with a hermitian metric $h'$, which we
want to compare with $h$, the hermitian metric on $V$ induced from the
monad (\ref{monad}). The key point is a fact we noted before in lemma
\ref{model}: given an 1-form $a$ on $\tproj$, its $L^2$-norm with respect
to the round metric is always larger than its $L^2$-norm with respect
to the flat metric on $T\times(\proj\setminus\{\infty\})$:
$$ ||a||_{L^2_R}>||a||_{L^2_F} $$
Thus, comparing the monads (\ref{monad}) and (\ref{cpx1}), one sees that
$h$ is bounded above by $h'$. In particular, the metric $h$ is bounded
at $\as$.

We can regard $\vv$ as an {\em index bundle} for the family of
Dirac operators over $\tproj$ parametrised by $\ksi\in\dual$.
Hence, its degree can be computed by the Atiyah-Singer index
theorem for families. Consider now the bundle
${\bf G}=p_{12}^*\ee\otimes p_{13}^*{\bf P}$ over $\tproj\times\dual$,
and note that ${\bf G}|_{\tproj\times\{\ksi\}}=\ee(\ksi)$. Then we
have:
\begin{eqnarray*}
ch(\vv) & = & - ch({\bf G})\cdot td(\tproj)/[\tproj] \ = \\
& = & - \left( 2+2c_1({\bf P})+c_1({\bf P})^2-c_2(\ee) \right)
\left( 1+\frac{1}{2}c_1(\proj) \right)/[\tproj] \ = \\
& = & k-\frac{1}{2}c_1({\bf P})^2c_1(\proj)/[\tproj] \ = \
k-2\hat{t}
\end{eqnarray*}
where the ``$-$" sign in the first line is needed since $\vv$ is
formed by the null spaces of the adjoint Dirac operator.

Summing up:
\begin{lem} \label{degree}
The dual bundle $(V,\del_B)\seta\dual\setminus\{\as\}$ admits a
holomorphic extension $\vv\seta\dual$ of degree $-2$. Moreover,
its hermitian metric $h$ is bounded above at the punctures $\as$.
\end{lem}

The determinant line bundle of $\vv$ is not fixed, however. In fact,
let $t_x:\tproj\seta\tproj$ be the translation of the torus by $x\in T$,
acting trivially on $\proj$, and let $\ee'=t_x^*\ee$. If $\vv'$ is
the dual bundle associated with $\ee'$ then $\vv'=\vv\otimes L_x$.
Indeed:
\begin{eqnarray*}
\vv_\ksi' = H^1(\tproj,\ee'(\ksi)) & = &
H^1\left( \tproj,p_{12}^*(t_x^*\ee)\otimes p_{13}^*{\bf P}|_{\tproj\times \{ \ksi \} } \right) \\
& = & H^1 \left( \tproj,t_x^*(p_{12}^*\ee\otimes p_{13}^*{\bf P})\otimes p_3^*L_x
      |_{\tproj\times \{ \ksi \} } \right) \\
& = & H^1 \left( \tproj,p_{12}^*\ee\otimes p_{13}^*{\bf P}|_{\tproj\times \{ \ksi \} } \right)
      \otimes (L_x)_\ksi \\
\Rightarrow \ \ \vv_\ksi' & = & \vv_\ksi\otimes(L_x)_\ksi
\end{eqnarray*}
as a canonical isomorphism for each $\ksi\in\dual$. Thus
$\vv'=\vv\otimes L_x$.

Note also that if $B$ is an admissible connection, $\vv$ admits no
splitting $\vv=\vv_0\oplus L$ compatible with $B$ for any flat line
bundle $L$.

%-------------------------------------------------------------------------

\paragraph{Defining the Higgs field.}
The next step is to give a holomorphic description of the Higgs field $\Phi$.

Recall that $h^0(\tproj,p_2^*\oo_{\proj}(1))=2$, and
regarding $\proj=\cpx\cup\{\infty\}$, we can fix two holomorphic
sections $s_0,s_\infty\in H^0(\proj,\oo_{\proj}(1))$ such that
$s_0$ vanishes at $0\in\cpx$ and $s_\infty$ vanishes at the point
added at infinity. In homogeneous coordinates
$\{(w_1,w_2)\in\cpx^2|w_2\neq0\}$ and $\{(w_1,w_2)\in\cpx^2|w_1\neq0\}$,
we have that, respectively ($w=w_1/w_2$):
\begin{eqnarray*}
s_0(w) = w & \ \ \ & s_0(w) = 1 \\
s_\infty(w) = 1 & \ \ \ & s_\infty(w) = \frac{1}{w}
\end{eqnarray*}

Let us first consider an alternative definition of the transformed Higgs
field. For each $\ksi\in\dual$, we define the map:
\begin{equation} \begin{array}{rcl}
H^1(\tproj,\ee(\ksi))\times H^1(\tproj,\ee(\ksi)) &
\stackrel{\Psi_\ksi}{\longrightarrow} & H^1(\tproj,\tilde{\ee}(\ksi)) \\
(\alpha,\beta) & \mapsto & \alpha\otimes s_0-\beta\otimes s_\infty
\end{array} \end{equation}
If $(\alpha,\beta)\in{\rm ker}\Psi_\ksi$, we define an endomorphism
$\varphi$ of $H^1(\tproj,\ee(\ksi))$ at the point $\ksi\in\dual$ as
follows:
\begin{equation} \label{alt.hig2}
\varphi_\ksi(\alpha)=\beta
\end{equation}

We check that $\varphi$ actually coincides with the
Higgs field $\Phi$ we defined in the previous section, up to a 
multiplicative constant. Note that:
$$ \alpha\otimes s_0-\beta\otimes s_\infty=0 \ \ \Leftrightarrow \ \
   \beta = \alpha(\otimes s_0)(\otimes s_\infty)^{-1} $$
Moreover, recall that, for any trivialisation of $\oo_{\proj}(1)$ with
local coordinate $w$ on $\proj$, the quotient $\frac{s_0(w)}{s_\infty(w)}=w$.
The claim now follows from the proof of proposition \ref{spin/coho};
we denote $\Phi_\ksi=2\sqrt{2}\pi\cdot \varphi_\ksi$.

\begin{prop} \label{higgs}
The eigenvalues of the Higgs field $\Phi$ have at most simple poles
at $\as$. Moreover, the residues of $\Phi$ are semi-simple and have
rank $\leq2$  if $\ksi_0$ is an element of order 2 in the Jacobian
of $T$, and rank $\leq1$ otherwise.
\end{prop}

\pf Suppose $\alpha(\ksi)$ is an eigenvector of $\Phi_\ksi$ with
eigenvalue $\epsilon'(\ksi)=1/\epsilon(\ksi)$, i.e.
$\Phi_\ksi(\alpha(\ksi)) = \epsilon'(\ksi)\cdot \alpha(\ksi)$.
Thus,
$$ \alpha(\ksi)\otimes s_0 - \epsilon'(\ksi)\cdot \alpha(\ksi)\otimes s_\infty =0
   \ \ \Rightarrow \ \
   \alpha(\ksi)\otimes (\epsilon(\ksi)\cdot s_0-s_\infty) =0$$
Therefore, denoting $s_\epsilon(\ksi)=\epsilon(\ksi)\cdot s_0-s_\infty$,
we have that $\alpha(\ksi)\in{\rm ker}(\otimes s_\epsilon(\ksi))$.

On the other hand, consider the sheaf sequence:
$$ 0 \seta \ee(\ksi) \stackrel{\otimes s_\epsilon(\ksi)}{\seta}
     \widetilde{\ee}(\ksi) \seta \widetilde{\ee}(\ksi)|_{T_{\epsilon'(\ksi)}} \seta 0 $$
since the section $s_\epsilon(\ksi)$ vanishes at $\epsilon'(\ksi)$. It
induces the cohomology sequence:
\begin{equation} \label{sqc2}
0 \seta H^0(T_{\epsilon'(\ksi)},\tilde{\ee}(\ksi)|_{T_{\epsilon'(\ksi)}})
  \seta H^1(\tproj,\ee(\ksi)) \stackrel{\otimes s_{\epsilon}(\ksi)}{\seta} ...
\end{equation}
so that ${\rm ker}(\otimes s_\epsilon(\ksi))=
H^0(T_{\epsilon'(\ksi)},\tilde{\ee}(\ksi)|_{T_{\epsilon'(\ksi)}})$
which is non-empty if and only if
$\ee(\ksi)|_{T_{\epsilon'(\ksi)}}=L_\ksi\oplus L_{-\ksi}$ or ${\bf F}_2\otimes L_\ksi$.

Hence, as $\ksi$ approaches $\as$, we must have that one of the eigenvalues
of $\Phi$, say $\epsilon'(\xi)$ approaches $\infty$,
since $\ee|_{T_\infty}=L_{\xi_0}\oplus L_{-\xi_0}$. Moreover,
$s_\epsilon(\ksi)\seta s_\infty$, so that:
$$ \lim_{\ksi\seta\as}\alpha(\ksi) \in {\rm ker}(\otimes s_\infty)=
   H^0(T_\infty,\ee(\ksi)|_{T_\infty}) $$

Therefore, we conclude that, if $\xi_0\neq-\xi_0$, then one of the eigenvalues of
$\Phi$ has a simple pole at $\as$ since $h^0(T_\infty,\ee(\as)|_{T_\infty})=1$;
similarly, if $\xi_0=-\xi_0$, then two of the eigenvalues of
$\Phi$ have a simple poles at $\ksi_0$.

Note in particular that the images of the residues
of $\Phi$ at $\as$ are precisely given by:
$$ H^0(T_\infty,\tilde{\ee}(\as)|_{T_\infty}) \subset H^1(\tproj,\ee(\as)) $$
\pfend

This proposition almost concludes the main task of this paper, namely
to construct the inverse of the Nahm transform of \cite{J2}. It only
remains to be shown that the Nahm transformed Higgs pair is admissible.
We must then show how to match the $SU(2)$ bundle $\check{E}\seta\torus$
with doubly-periodic instanton $\check{A}$ constructed from the transformed
Higgs pair $(B,\Phi)$ as in \cite{J2} with the original objects $A$ and
$E\seta\torus$ we started with in the present paper. These tasks are
taken up in the following section.

%--------------------------------------------------------------------

\section{Proof of inversion} \label{inv}

So far, we have established that the Nahm transform of a doubly-periodic
instantons is the same kind of singular Higgs pair as those we started
with in the first part of this series \cite{J2}.

We must now show that the transform presented here is actually the
inverse of the construction of instantons of \cite{J2}. More precisely,
we show that if we start with a doubly-periodic instanton $A$, apply
the Nahm transform to obtain a Higgs pair $(B,\Phi)$, then the
corresponding doubly-periodic instanton constructed as in \cite{J2} is
gauge equivalent to the original object.

First, consider the six-dimensional manifold $\torus\times(\dual\setminus\{\as\})$.
To shorten notation, we denote $M_\ksi=\torus\times\{\ksi\}$ and
$\dual_{(z,w)}=\{z\}\times\{w\}\times(\dual\setminus\{\as\})$.

Now take the bundle $\calg=p_{23}^*E\otimes p_{12}^*\mathbf{P}$
over $\torus\times(\dual\setminus\{\as\})$; note that
$\calg|_{M_\ksi}\equiv E(\ksi)$ and
$\calg|_{\dual_{(z,w)}} \equiv \underline{E_{(z,w)}}\otimes L_z$,
where $\underline{E_{(z,w)}}$ denotes a trivial rank 2 bundle over
$\dual\setminus\{\as\}$ with the fibres canonically identified with
the vector space $E_{(z,w)}$.

$\cal G$ is clearly holomorphic; we denote by $\del_M$ the action
of the associated Dolbeault operator along the $\tproj$ direction,
and by $\del_{\dual}$ its action along the $\dual$ direction. In
particular, $\del_M|_{M_\ksi} \equiv \del_{A_\ksi}$.

Let ${\bf C}^{p,q}=\Lambda^{0,p}_{\torus}({\cal G})\otimes\Lambda^q_{\dual}({\cal G})$;
in other words, ${\bf C}^{p,q}$ consists of the $(p+q)$-forms over
$\torus\times(\dual\setminus\{\as\})$ with values in $\cal G$ spanned by
forms of the shape:
\begin{equation} \label{cpq} \begin{array}{c}
s(z,w,\ksi)d\overline{z}_{i_1}d\overline{w}_{i_2}d\ksi_{j_1}d\overline{\ksi}_{j_2}, \\
i_1,i_2,j_1,j_2\in\{0,1\} \ \ \ {\rm and} \ \ \ i_1+i_2=p,\ j_1+j_2=q
\end{array} \end{equation}
Analytically, we want to regard ${\bf C}^{p,q}$ as the completion of the
set of smooth forms of the shape above with respect to a Sobolev norm
described as follows:
$$ \begin{array}{cl}
\left| s|_{\torus\times\{\ksi\}} \right| \in L^2_q(\Lambda^{2-q}E(\ksi)) &
{\rm for\ each}\ \ksi\in\dual\setminus\{\as\} \\
\left| s|_{\{(z,w)\}\times\dual\setminus\{\as\}} \right| \in L^2_q(\Lambda^{2-q}L_z) &
{\rm for\ each}\ (z,w)\in\torus
\end{array} $$

Now, define the maps:
\begin{equation} \label{this} \begin{array}{rcl}
\mathbf{C}^{p,0} \stackrel{\delta_1}{\seta} & {\bf C}^{p,1} &
\stackrel{\delta_2}{\seta} \mathbf{C}^{p,2} \\
\delta_1(s)=(\del_{\dual}s,-w\cdot s\wedge d\ksi) & &
\delta_2(s_1,s_2)=(\del_{\dual}s_2+w\cdot s_1\wedge d\ksi)
\end{array} \end{equation}
for $(s_1,s_2)\in\Lambda^{0,p}_{\tproj}({\cal G})\otimes
\left( \Lambda^{0,1}_{\dual}({\cal G})\oplus\Lambda^{1,0}_{\dual}({\cal G}) \right)
\equiv {\bf C}(p,1)$. Note that (\ref{this}) does define a complex.

The inversion result will follow from the analysis of the spectral
sequences associated to the following double complex
(for the general theory of spectral sequences and double complexes,
we refer to \cite{BT}):
\begin{equation} \label{dblcpx} \begin{array}{ccccc}
\mathbf{C}^{0,2} & \stackrel{\del_M}{\seta} & \mathbf{C}^{1,2} &
\stackrel{-\del_M}{\seta} & \mathbf{C}^{2,2} \\
\ \ \uparrow \delta_2 & & \ \ \ \uparrow -\delta_2 & & \ \ \uparrow \delta_2 \\
\mathbf{C}^{0,1}& \stackrel{\del_M}{\seta} &
\mathbf{C}^{1,1}& \stackrel{-\del_M}{\seta} & \mathbf{C}^{2,0} \\
\ \ \uparrow \delta_1 & & \ \ \ \uparrow -\delta_1 & & \ \ \uparrow \delta_1 \\
\mathbf{C}^{0,0} & \stackrel{\del_M}{\seta} & \mathbf{C}^{1,0} &
\stackrel{-\del_M}{\seta} & \mathbf{C}^{2,0}
\end{array} \end{equation}
The idea is to compute the total cohomology of the spectral sequence
in the two possible different ways and compare the filtrations of
the total cohomology.

\begin{lem} \label{rows}
By first taking the cohomology of the rows, we obtain
\begin{equation} \begin{array}{cccccc}
          & \ & & 0 & H^2({\cal C}(e,0)) & 0 \\
E^{p,q}_2 & \ & & 0 & H^1({\cal C}(e,0)) & 0 \\
& \ & q\uparrow & 0 & H^0({\cal C}(e,0)) & 0 \\
& \ & & \stackrel{\seta}{p} & &
\end{array} \end{equation}
where $H^i({\cal C}(e,0))$ are the cohomology groups of the complex
that yields the monad description of the construction of doubly-periodic
instantons in \cite{J2} (see proposition 3 there).
\end{lem}

\pf First, note that the rows coincide with the complex (\ref{monad}).

Moreover, we can regard elements in $\mathbf{C}^{p,q}$ as $q$-forms over
$\dual$ with values in $L^2_{2-p}(\Lambda^{0,p}_{\torus}{\cal G})$.
To see this, fix some $\ksi'\in\dual$; by (\ref{cpq}),
$s(z,w,\ksi')\in\Lambda^{0,p}{\cal G}|_{M_{\ksi'}}$. So,
by varying $\ksi'$ we get the interpretation above.

This said, it is clear that the first and second columns of
$E^{p,q}_1$ must vanish, since $A$ is irreducible. In the middle column,
we get $q$-forms over $\dual$ with values in ${\rm ker}(\del_M^*-\del_M)$,
which for a fixed $\ksi'$ restricts to ${\rm ker}(D^*_{A_{\ksi'}})$.

Therefore, after taking the cohomologies of the rows, we are left with:
\begin{equation} \label{rows2} \begin{array}{cccccc}
                & \ & & 0 \ \ & L^p(\Lambda^{1,1}V) & \ \ 0 \\
                & \ & &  & \ \ \ \ \ \ \uparrow(\del_B+\Phi) &  \\
{\bf C}^{p,q}_1 & \ & & 0 \ \ & L^p_1(\Lambda^{1,0}V\oplus\Lambda^{0,1}V) &\ \ 0 \\
                & \ & &  & \ \ \ \ \ \ \uparrow(\Phi+\del_B) &  \\
                & \ & q\uparrow & 0 \ \ & L^p_2(\Lambda^0V) &\ \ 0 \\
                & \ & & \stackrel{\seta}{p} & &
\end{array} \end{equation}
But this is just the complex that yields the monad description of
the construction of doubly-periodic instantons in \cite{J2}. The
lemma follows after taking the cohomology of the remaining column.
\pfend

\paragraph{Total cohomology and admissibility.}
Note that, as we pointed out in the beginning of this section, we
still do not know if the Higgs pair $(B,\Phi)$ arising from the
instanton $(E,A)$ is admissible or not, i.e. the hypercohomology
spaces ${\Bbb H}^0$ and ${\Bbb H}^2$ might be nontrivial. The next
lemma deals with this problem.

\begin{lem} \label{t.coh}
The only nontrivial cohomology of the total complex is \linebreak
$H^2({\bf C}(p,q))$, which is naturally isomorphic to the fibre
$E_{(e,0)}$.
\end{lem}

In particular, this shows that the Higgs pairs $(B,\Phi)$ obtained
via Nahm transform on instanton connection $A\in\cfg$ are indeed
admissible, see \cite{J2}.

\pf First note that we can regard an element in ${\bf C}^{p,q}$ as a
$(0,p)$-form over $\torus$ with values in $\Lambda^{q_1,q_2}_{\dual}(\calg)$.
Since $\calg|_{\dual_{(z,w)}} \equiv \underline{E_{(z,w)}}\otimes L_z$,
${\rm ker}\del_M$ and ${\rm ker}\del_M^*$ are nontrivial only if $z=e$,
the identity element in the group law of $T$. Hence, it is enough to work on a
tubular neighbourhood of $\{e\}\times\proj(\dual\setminus\{\as\})$.

More precisely, we define another double complex $({\rm germ}\ {\bf C})^{p,q}$,
consisting of forms defined on arbitrary neighbourhoods of
$\{e\}\times\proj\times(\dual\setminus\{\as\})$. Then we have a restriction
map ${\bf C}^{p,q}\seta({\rm germ}\ {\bf C})^{p,q}$ commuting with
$\del_M$, $\delta_1$ and $\delta_2$. Such map also induces an
isomorphism between the total cohomologies of ${\bf C}^{p,q}$ and
$({\rm germ}\ {\bf C})^{p,q}$. So we can work with $({\rm germ}\ {\bf C})^{p,q}$
to prove the lemma.

Let $V_e$ be some neighbourhood of $e\in T$. By the Poincar\'e lemma
applied to $\del_T$, we get:
\begin{equation} \label{t.coh2} \begin{array}{ccccccc}
& \ & & & \Lambda^2_{V_e}(\calg) & 0 & 0 \\
& \ & &  & \ \ \ \ \ \ \uparrow & & \\
({\rm germ}\ {\bf C})^{p,q}_1
& \ & & & \Lambda^1_{V_e}(\calg) & 0 & 0 \\
& \ & &  & \ \ \ \ \ \ \uparrow & &  \\
& \ & q\uparrow & & \Lambda^0_{V_e}(\calg) & 0 & 0 \\
& \ & & \stackrel{\seta}{p} & & &
\end{array} \end{equation}
where $V_e$ denotes a tubular neighbourhood of
$N_e=\{e\}\times\proj\times(\dual\setminus\{\as\})$

As in \cite{DK} (see pages 91-92), the complex in the first row is,
after restriction, mapped into a Koszul complex over $N_e$:
$$ \oo_{N_e}(\calg) \stackrel{(w\ \ksi)}{\longrightarrow}
   \oo_{N_e}(\calg)\oplus \oo_{N_e}(\calg)
   \stackrel{(-\ksi,z)}{\seta}\oo_{N_e}(\calg) $$
so that:
\begin{equation} \label{t.coh3} \begin{array}{cccccc}
& \ & & E_{(e,0)} & 0 & 0 \\
({\rm germ}\ {\bf C})^{p,q}_2
& \ & & 0 & 0 & 0 \\
& \ & q\uparrow & 0 & 0 & 0 \\
& \ & & \stackrel{\seta}{p} & &
\end{array} \end{equation}
\pfend

It then follows from lemmas \ref{rows} and \ref{t.coh} that there is a
natural isomorphism of vector spaces
${\cal I}_I:H^1({\cal C}(e,0))\equiv\check{E}_{(e,0)}\seta E_{(e,0)}$,
which in principle may depend on the choice of complex structure $I$ on $\torus$.

\paragraph{Matching $\mathbf{(\check{E},\check{A})}$ with the original data.}
Since the choice of identity element in $T$ and of origin in $\cpx$
is arbitrary, we can extend ${\cal I}_I$ to a bundle isomorphism
$E\seta\check{E}$. More precisely, let $t_{(u,v)}:\torus\seta\torus$
be the translation map $(z,w)\seta(z+u,w+v)$. Clearly, the connection
$t_{(u,v)}^*A$ on the pullback bundle $t_{(u,v)}^*E$ is also irreducible
and $t_{(u,v)}^*E_{(e,0)}\equiv E_{(u,v)}$. Computing the total
cohomology of the double complex (\ref{dblcpx}) associated to the
bundle $t_{(u,v)}^*{\cal G}$ (where $t_{(u,v)}^*$ acts trivially on
$\dual$ coordinate), lemmas \ref{rows} and \ref{t.coh} lead to an isomorphism
of vector spaces $H^1({\cal C}(u,v))\equiv\check{E}_{(u,v)}\seta E_{(u,v)}$.

It is clear from the naturality of the constructions that these
fibre isomorphisms fit together to define a holomorphic bundle
isomorphism \linebreak ${\cal I}_I:E\seta\check{E}$. In particular,
${\cal I}_I$ takes the Dolbeault operator $\del_A$ of the
holomorphic bundle $E\seta\torus$ to the Dolbeault operator
$\del_{\check{A}}$ of the holomorphic bundle
$\check{E}\seta\torus$. It also follows from this observation that
the holomorphic extensions $\ee$ and $\check{\ee}$ must be
isomorphic as holomorphic vector bundles.

However, such fact still does not guarantee that the connections
$A$ and $\check{A}$ are gauge-equivalent. This is accomplished if
we can show that ${\cal I}_I$ is actually {\em independent} of the
choice of complex structure in $\torus$. Therefore, the proof of
the main theorem \ref{nahmthm} is completed by the following
proposition:

\begin{prop}
The bundle map ${\cal I}_I:\check{E}\seta E$ is independent of the
choice of complex structure on $\torus$.
\end{prop}

\pf Again, it is sufficient to consider only the fibre over
$(e,0)$. As in \cite{DK}, the idea is to present an explicit
description of ${\cal I}_I:\check{E}_{(e,0)}\seta E_{(e,0)}$,
and then show that it is Euclidean invariant.

Let $\alpha\in H^1({\cal C}(e,0))\subset{\bf C}^{1,1}$. To find
${\cal I}_I([\alpha])$ we have to find $\beta\in{\bf C}^{0,2}$ such that
$\del_M\beta=\delta_2\alpha$. A solution to this equation is provided
by the Hodge theory for the $\del_M$ operator:
$$ \beta=G_M(\del_M^*\delta_2\alpha) $$
where $G_M$ denotes the Green's operator for $\del_M^*\del_M$,
which can be regarded fibrewise as the family of Green's operators
$G_{A_\ksi}=G_M|_{M_\ksi}$ parametrised by $\ksi\in(\dual\setminus\{\as\})$.

In principle, $\beta$ depends on the complex structure $I$ via the
operators $\del_M$ and $G_M$. However, by the Weitzenb\"ock formula
applied to the bundle $\calg$, we have:
$$ \del_M^*\del_M=\nabla^*_M\nabla_M $$
Here, $\nabla_M$ is the covariant derivative in the $\torus$ direction
on $\calg$. With this interpretation, $G_M=(\nabla^*_M\nabla_M)^{-1}$
is seen to be independent of the complex structure $I$; in fact, it is
Euclidean invariant.

Now $\beta$ as an element of ${\bf C}^{1,1}$ has the form
$\beta(z,w;\ksi)d\ksi d\overline{\ksi}$, so that the restriction
$r_{(e,0)}(\beta)=\beta|_{\dual_{(e,0)}}$ is a $(1,1)$-form over
$\dual\setminus\{\as\}$ with values in $E_{(e,0)}$. Take its
cohomology class in $H^2(\dual\setminus\{\as\},\underline{\cpx}\otimes\underline{E_{(e,0)}})$,
so that:
$$ {\cal I}_I([\alpha])=\int_{\dual_{(e,0)}}r_{(e,0)}(\beta) $$
which is the desired explicit description. \pfend

Summing up the work done in \cite{J2} and in this paper, we
have proven theorem \ref{nahmthm}.

%---------------------------------------------------------------

\section{Instantons of higher rank}

One easily realizes that there is nothing really
special about rank two bundles; the whole proof could easily be
generalised to higher rank. Indeed, the only point in restricting to the
rank two case is to reduce the number of possible vector bundles over an
elliptic curve, and avoid a tedious case-by-case study throughout the
various stages of the proof.

Before we can state the generalisation of the main theorem \ref{nahmthm},
we must review our definitions of asymptotic state and irreducibility.

The restriction of the holomorphic extension $\ee\seta\tproj$ to the
added divisor $T_\infty$ is a flat $SU(n)$ bundle, i.e.
$$ \overline{E}|_{T_\infty}=L_{\ksi_1}\oplus\dots\oplus L_{\ksi_k} $$
$$ {\rm such\ that}\ \ \ \ \ \bigotimes_{l=1}^k L_{\ksi_l}=\oo_{T} $$
In other words, $\overline{E}|_{T_\infty}$ is determined by a set
of points $(\ksi_1,\dots,\ksi_j)\in{\cal J}(T)$ with multiplicities
$(m_1,\dots,m_j)$, and such that $\sum_{l=1}^j m_l\ksi_l=0$. We
call such data the {\em generalised asymptotic state}.

Moreover, we will say that $(E,A)$ is {\em 1-irreducible} if there is
no flat line bundle $E\seta\torus$ such that $E$  admits a splitting
$E'\oplus L$ which is compatible with the connection $A$.

\begin{thm} \label{nahmthmgen}
There is a bijective correspondence between the following objects:
\begin{itemize}
\item gauge equivalence classes of 1-irreducible $SU(n)$-instantons over $\torus$
with fixed instanton number $k$ and generalised asymptotic state
\linebreak 
$(\ksi_1,\dots,\ksi_j)$ with multiplicities $(m_1,\dots,m_j)$;
\item admissible $U(k)$ solutions of the Hitchin's equations over the dual torus
$\dual$, such that the Higgs field has at most simple poles at
$\{\ksi_1,\dots,\ksi_j\}$; moreover, its residue at $\ksi_j$ is
semi-simple and has rank $\leq m_j$.
\end{itemize}
\end{thm}

Also, the same remark about the possibility of
removing the technical hypothesis on the non-triviality of the
asymptotic states holds.

%--------------------------------------------------------------------

\section{Conclusion}

The attentive reader might have noticed that the assumptions on 
doubly-periodic instantons used on this paper (namely extensibility)
do coincide with the conclusions of the first paper of the series.
However, it is important to point out at this stage the small gap 
remaining between the conclusions of the present paper and the 
assumptions in \cite{J2}. 

More precisely, we assumed in \cite{J2} that the harmonic metric 
associated with the Higgs pair $(B,\Phi)$ on the bundle 
$V\seta\dual\setminus\{\as\}$ is non-degenerate along the kernel 
of the residues of $\Phi$, and $h\sim O(r^{1\pm\alpha})$ along 
the image of the residues of $\Phi$, for some alpha 
$0\leq\alpha<1/2$, in a holomorphic trivialisation 
of $V$ over a sufficiently small neighbourhood around $\as$, . 

The gap is closed in \cite{BJ}, where it is shown that the Nahm 
transformed Higgs pairs here constructed do satisfy the above 
condition.

The analytical features of extensible doubly-periodic instantons are 
further studied by Olivier Biquard and the author in \cite{BJ}. In particular, 
we have provided a deformation theory description of the moduli space of 
rank two doubly-periodic instanton connections as a hyperk\"ahler manifold 
of dimension $4k-2$. Moreover, it is also shown that the Nahm transform is
a hyperk\"ahler isometry between the moduli of doubly-periodic instantons 
and the moduli of singular Higgs pairs.

%--------------------------------------------------------------------

 \end{document}